\documentclass{ifacconf}

\usepackage{amssymb,amsfonts}
\usepackage{amsmath}
\usepackage{graphicx}      
\usepackage{natbib}        
\usepackage{mathtools}
\usepackage{xcolor}
\usepackage{bm}
\usepackage{comment}
\usepackage{tikz}
\usepackage{ifthen}
\usepackage{algorithm}
\usepackage{algpseudocode}
\usepackage{url}
\usepackage{float}
\usepackage{listings}
\newfloat{algorithm}{t}{lop}
\usepackage{balance}

\newcommand{\defeq}{:=}

\let\lamon\l
\renewcommand{\l}{_\textrm{lin}}

\newcommand{\Z}{\mathbf{z}}
\newcommand{\V}{\mathbf{v}}
\newcommand{\U}{\mathbf{u}}
\newcommand{\M}{M}

\renewcommand{\P}{\mathbf \Phi}

\newcommand{\X}{\mathbf{x}}

\newcommand{\nx}{n_\mathrm{x}}
\renewcommand{\nu}{n_\mathrm{u}}
\newcommand{\nc}{n_\mathrm{c}}
\newcommand{\nf}{n_\mathrm{f}}

\newcommand{\nw}{n_\mathrm{w}}
\newcommand{\ny}{{n_\mathrm{y}}}
\newcommand{\R}{\mathbb{R}}

\newcommand{\Qf}{P}

\newcommand{\f}{\mathrm{f}}

\newcommand{\T}{^\top}

\newcommand{\x}{\mathrm{x}}
\renewcommand{\u}{\mathrm{u}}

\newcommand{\Px}{\bm{\Phi}^{\mathrm{x}}}
\newcommand{\Pu}{\bm{\Phi}^{\mathrm{u}}}
\newcommand{\F}{\mathcal{F}}

\newcommand{\bo}{^\mathrm{ct}}

\newif\ifextended
\extendedtrue
\ifextended
\newcommand{\extension}[2]{#1}
\else
\newcommand{\extension}[2]{#2}
\fi

\newtheorem{remark}{Remark}[section]

\usepackage{glossaries}
\newacronym{MPC}{MPC}{model predictive control}
\newacronym{SLS}{SLS}{system level synthesis}
\newacronym{LTV}{LTV}{linear time-varying}
\newacronym{NLP}{NLP}{nonlinear program}
\newacronym{SQP}{SQP}{sequential quadratic programming}
\newacronym{SCP}{SCP}{sequential convex programming}

\newacronym{SOCP}{SOCP}{second-order cone program}
\newacronym{QP}{QP}{quadratic program}
\newacronym{LQR}{LQR}{linear quadratic regulator}
\newacronym{KKT}{KKT}{Karush–Kuhn–Tucker}

\begin{document}
\begin{frontmatter}

\title{Fast System Level Synthesis:\\ Robust Model Predictive Control using Riccati Recursions\thanksref{footnoteinfo}} 

\thanks[footnoteinfo]{
This work has been supported by the European Space Agency under OSIP 4000133352, the Swiss Space Center, the Swiss National Science Foundation under NCCR Automation (grant agreement 51NF40 180545), 
by DFG via project 424107692, and
by the European Union's Horizon 2020 research and innovation programme, Marie~Sk\lamon{}odowska-Curie grant agreement No. 953348,~\mbox{ELO-X}. 
}

\author[First]{Antoine P. Leeman}
\author[First]{Johannes K{\"o}hler}
\author[Second]{Florian Messerer} %
\author[First]{Amon Lahr}
\author[Second,Third]{Moritz Diehl}
\author[First]{Melanie N. Zeilinger}

\address[First]{Institute for Dynamic Systems and Control, ETH Zurich, Zurich 8053, Switzerland (e-mail:\{aleeman; jkoehle; amlahr; mzeilinger\}@ethz.ch)}
\address[Second]{Department of Microsystems Engineering (IMTEK), University
of Freiburg, 79110 Freiburg, Germany (e-mail: \{florian.messerer; moritz.diehl\}@imtek.uni-freiburg.de)}
\address[Third]{Department of Mathematics, University of Freiburg, 79104 Freiburg, Germany}

\begin{abstract}              
System level synthesis enables improved robust MPC formulations by allowing for joint optimization of the nominal trajectory and controller. This paper introduces a tailored algorithm for solving the corresponding disturbance feedback optimization problem for linear time-varying systems. The proposed algorithm iterates between optimizing the controller and the nominal trajectory while converging q-linearly to an optimal solution.  We show that the controller optimization can be solved through Riccati recursions leading to a horizon-length, state, and input scalability of $\mathcal{O}(N^2 ( \nx^3 +\nu^3))$ for each iterate. On a numerical example, the proposed algorithm exhibits computational speedups by a factor of up to $10^3$ compared to general-purpose commercial solvers.
\end{abstract}
\glsreset{SLS}
\glsreset{MPC}

\end{frontmatter}

\section{Introduction}

The ability to plan safe trajectories in real-time, despite model mismatch and external disturbances, is a key enabler for high-performance control of autonomous systems~\citep{majumdar2017funnel,schulman2014motion}.
{iLQR}~\citep{li2004iterative,giftthaler2018family} methods have shown great practical benefits~\citep{neunert2016fast, howell2019altro,singh2022optimizing} by producing a feedback controller based on the optimal nominal trajectory.
However, robust stability and robust constraint satisfaction are generally only addressed heuristically~\citep{Manchester2017}.

\Gls{MPC}~\citep{rawlings2017model} has emerged as a key technology for control subject to (safety) constraints. 
Real-time feasibility of \gls{MPC}, even in fast-sampled systems, is largely due to advancements in numerical optimization~\citep{frison2016algorithms,domahidi2012efficient}. These tailored algorithms exploit the stage-wise structure of the nominal trajectory optimization problem, utilizing, e.g., efficient Riccati recursions.
Robust \gls{MPC} approaches can account for model mismatch and external disturbances by propagating them over the prediction horizon.
While few robust \gls{MPC} approaches optimize a feedback policy~\citep{Scokaert1998,villanueva2017robust,messerer2021efficient,kim2022joint}, most robust MPC designs use an offline-\textit{fixed} controller~\citep{mayne2005robust,Kohler2021ASystems,zanelli2021zero} for computational efficiency.

Disturbance feedback for \gls{LTV} systems \citep{goulart2006optimization,ben2004adjustable} overcomes this limitation by using a convex controller parametrization.
Building on this parametrization, \gls{SLS}~\citep{anderson2019system,tseng2020system,li2023learning} facilitates the consideration of (structural) constraints on the closed-loop response.
\Gls{SLS} has also been proposed within an \gls{MPC} formulation~\citep{sieber2021system, chen2022robust,Leeman2023RobustSynthesis,leeman2023_CDC}.
However, the large number of decision variables in \gls{SLS} poses a major obstacle for real-time implementation.
\begin{figure}[ht!]
    \centering
\tikzset{every picture/.style={line width=0.75pt}} 
\begin{tikzpicture}[x=0.75pt,y=0.75pt,yscale=-1,xscale=1]

\draw    (80,40) -- (80,20) -- (230,20) -- (230,38) ;
\draw [shift={(230,40)}, rotate = 270] [color={rgb, 255:red, 0; green, 0; blue, 0 }  ][line width=0.75]    (10.93,-3.29) .. controls (6.95,-1.4) and (3.31,-0.3) .. (0,0) .. controls (3.31,0.3) and (6.95,1.4) .. (10.93,3.29)   ;
\draw    (80,92) -- (80,110) -- (230,110) -- (230,90) ;
\draw [shift={(80,90)}, rotate = 90] [color={rgb, 255:red, 0; green, 0; blue, 0 }  ][line width=0.75]    (10.93,-3.29) .. controls (6.95,-1.4) and (3.31,-0.3) .. (0,0) .. controls (3.31,0.3) and (6.95,1.4) .. (10.93,3.29)   ;
\draw   (20,50) .. controls (20,44.48) and (24.48,40) .. (30,40) -- (140,40) .. controls (145.52,40) and (150,44.48) .. (150,50) -- (150,80) .. controls (150,85.52) and (145.52,90) .. (140,90) -- (30,90) .. controls (24.48,90) and (20,85.52) .. (20,80) -- cycle ;
\draw   (160,50) .. controls (160,44.48) and (164.48,40) .. (170,40) -- (280,40) .. controls (285.52,40) and (290,44.48) .. (290,50) -- (290,80) .. controls (290,85.52) and (285.52,90) .. (280,90) -- (170,90) .. controls (164.48,90) and (160,85.52) .. (160,80) -- cycle ;

\draw (85,65) node   [align=left] {\begin{minipage}[lt]{88.4pt}\setlength\topsep{0pt}
\begin{center}
{\scriptsize \textbf{Optimize}}\\{\scriptsize nominal trajectory $\displaystyle (\mathbf{z,v})$}\\{\scriptsize using QP}
\end{center}

\end{minipage}};
\draw (225,65) node   [align=left] {\begin{minipage}[lt]{88.4pt}\setlength\topsep{0pt}
\begin{center}
{\scriptsize \textbf{Optimize}}\\{\scriptsize controller gains $\displaystyle K_{k,j}$ }\\{\scriptsize via Riccati recursions}
\end{center}

\end{minipage}};
\draw (155,10) node  [font=\scriptsize] [align=left] {\begin{minipage}[lt]{129.2pt}\setlength\topsep{0pt}
\begin{center}
Update cost via duals $\displaystyle \mu $
\end{center}

\end{minipage}};
\draw (155,120) node  [font=\scriptsize] [align=left] {\begin{minipage}[lt]{129.2pt}\setlength\topsep{0pt}
\begin{center}
Update constraint tightening $\displaystyle \beta $
\end{center}

\end{minipage}};

\end{tikzpicture}

    \caption{
    Illustration how the robust MPC problem~\eqref{eq:sls} is solved by alternating between a nominal trajectory optimization problem with tightened constraints and the computation of optimal feedback gains $K_{k,j}$ using Riccati recursions.
    }
    \label{fig:iterative_update_sls}
\end{figure}
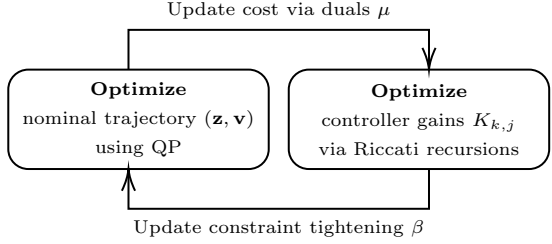
Asynchronous updates~\citep{sieber2023asynchronous} or {GPU} parallelization~\citep{alonso2022effective} can reduce the computational time.
\cite{goulart2008efficient} derived a structure-exploiting solver for disturbance feedback optimization, reducing the scalability of a naive interior-point method implementation from $\mathcal{O}(N^4)$ to $\mathcal{O}(N^3)$ per iteration.
In this work, we leverage results in numerical optimization~\citep{WrightStephenandNocedal1999NumericalOptimization,messerer2021efficient} to propose a customized optimization solver for \gls{SLS} with $\mathcal{O}(N^2)$ scalability, significantly improving the computation times.
\subsubsection*{Contribution}
This work addresses robust MPC for \gls{LTV} systems using an \gls{SLS} formulation (Section \ref{sec:setup}).
Our key contribution is an efficient solution achieved through iterative optimization of the nominal trajectory and the controller (Section \ref{sec:riccati_linear_sls}), as inspired by~\cite{messerer2021efficient} and illustrated in Fig.~\ref{fig:iterative_update_sls}.

The proposed algorithm has the following properties:
\begin{itemize}
    \item It converges q-linearly to an optimal solution of a tightened robust MPC problem.
\item Each controller optimization is solved with efficient Riccati recursions with a combined scalability of $\mathcal{O}(N^2 ( \nx^3 + \nu ^3 ) )$.
\item The proposed algorithm also allows for parallelization, leading to a scalability of $\mathcal{O}(N ( \nx^3 + \nu ^3 ) )$.
\end{itemize}

Furthermore, we discuss how the approach can be extended to nonlinear systems via \gls{SCP} (Remark~\ref{rmq:nonlinear_extension}).
We showcase the numerical properties of our solver for a range of state dimensions and horizon lengths in comparison to the commercial general-purpose solvers \cite{gurobi}, and Mosek~\cite{mosek} (Section~\ref{sec:numerical}).
\subsubsection*{Notation}
We denote stacked vectors {or matrices} by $(a,\ldots,b) = [a\T,\ldots,b\T]\T$. 
The matrix ${I}$ denotes the identity with its dimensions either inferred from the context or indicated by the subscript, i.e., ${I}_{\nx}\in\R^{\nx\times \nx}$, the matrix $\mathbf{0}_{n,m}\in \R^{n\times m}$ has all its elements zeros, and the vector $\bm{1}_{n}\in\R^{n}$ has all its elements one.
For a sequence of matrices $\Pu_{k,j}$, we define $\Pu_{0 :k, j}\defeq (\Pu_{0, j},\ldots, \Pu_{k, j})$.
Similarly, for a sequence of vectors $x_k$, we define $\mathbf{x} \defeq (x_0,\ldots, x_k)$.
We denote the 2-norm, resp. induced 2-norm, of a vector, resp. of a matrix, as $\|\cdot\|_2$.
We define the Frobenius norm of a matrix $A\in\R^{m\times n}$, as $\|A\|_\mathcal{F} = \text{Trace}(A\T A)$.
For a non-negative vector $x\in\R^n$, the square root $\sqrt{x}$ is to be understood elementwise. Besides, $0\le x\perp y \le 0$ is the shorthand notation for the three conditions \mbox{$x_i\cdot y_i=0,$}$~ x_i\ge 0,~ y_i\le 0,~ i=1,\ldots,n$.
The set of positive definite matrices is denoted by $\mathbb{S}_{++}^n$.
We denote the unit ball, centered at the origin, $\mathcal{E}_n:=\left\{ x \in \mathbb{R}^n, \|x\|_2\le 1\right\}$. For a vector-valued function $f: \mathbb{R}^n \rightarrow \mathbb{R}^q$, we denote by $\nabla \phi(x)\in \R^{n\times q}$ the transposed Jacobian, i.e., $\nabla \phi(x)= \frac{\partial \phi(x)}{\partial x}\T$.
\vspace{1cm}
\section{Problem setup \& System Level Synthesis}
\label{sec:setup}
We consider the following uncertain \gls{LTV} system:
\begin{equation}
    x_{k+1}=A_k x_k+B_k u_k+ E_k w_k,~x_0=\bar{x}_0,\label{eq:LTV}
\end{equation}
with state $x_k \in \mathbb{R}^{\nx}$, input $u_k \in \mathbb{R}^{n_{\mathrm{u}}}$ and disturbance lying in a unit ball $w_k \in \mathcal{E}_{\nx} \defeq \{w\in \R^{\nx}, \|w\|_2 \le 1 \}$.
To simplify the exposition, we assume that the disturbance scaling matrix $E_k$ is invertible. Remark~\ref{rem:extension_matrix_E} discusses the extension to left-invertible matrices $E_k \in \R^{\nx \times n_\mathrm{w}}$ with $w_k \in \R^{\nw}$.
The initial condition is given by $\bar{x}_0 \in \mathbb{R}^{n_{\mathrm{x}}}$.

We consider the problem of designing an optimal controller minimizing the cost function
\begin{equation}
 J(\X,\U) =  \sum_{k=0}^{N-1} {x}_k^{\top} Q {x}_k+{u}_k^{\top} R {u}_k+{x}_N^{\top} \Qf {x}_N, 
\end{equation}
with $Q\in \mathbb{S}_{++}^{\nx}$, $R\in\mathbb{S}_{++}^{\nu}$, $P\in \mathbb{S}_{++}^{\nx}$while robustly satisfying the constraints
    \begin{subequations}
    \label{eq:constraints}
    \begin{align}
    g_{k,i}\T(x_k, u_k)+b_{k,i} \le 0,~\forall i=1,\ldots, \nc,\\
    ~g_{\f,i}\T x_N+b_{\f,i} \le 0,~\forall i=1,\ldots, \nf,
    \end{align}
    \end{subequations}
for $ k=0,\ldots,N-1$ with $g_{k,i}\in \R^{\nx + \nu}$, and $g_{\f,i}\in \R^{\nx}$.

We optimize over disturbance feedback
\begin{equation}
\label{eq:disturbance_feedback}
u_k = v_k + \sum_{j=0}^{k-1} \Pu_{k,j} w_j,
\end{equation}
with nominal input $v_k \in \R^{\nu}$, i.e.,
we assign a distinct disturbance feedback matrix $\Pu_{k,j}\in \R^{\nu\times \nx}$ for each disturbance $w_j$ and control $u_k$ with $k>j$.
As common in \gls{SLS}, the resulting state sequence can be expressed as
\begin{equation}
x_k= z_k   + \sum_{j=0}^{k-1} \Px_{k,j} w_j,~ z_0 = x_0,
\end{equation}
where $z_k$ is the corresponding nominal state and the corresponding closed-loop response \mbox{$\Px_{k,j}\in \R^{\nx\times \nx}$} denotes the influence of the disturbance $w_j$ on the state $x_k$.
Starting with $\Px_{j+1,j} = E_j$, the propagation of the disturbance, commonly used in \gls{SLS}, is governed by
\begin{equation}
\label{eq:tube_propagation}
\begin{aligned}
\Px_{k+1,j} = A_k \Px_{k,j} + B_k\Pu_{k,j},\\
\end{aligned}
\end{equation}
$j=0,\ldots,N-1$, $k=j+1,\ldots, N-1$.

Combining the disturbance propagation for each time step~\eqref{eq:tube_propagation}, we formulate the robust optimal control problem via \gls{SLS} as a \gls{SOCP},
\begin{subequations}
\label{eq:sls}
\begin{align}
\min_{\substack{{\Px},{\Pu},\\\Z, \V}} & J(\Z,\V) +  \tilde H_0(\P),\label{eq:SLO}\\
\text { s.t. } ~ & {z}_{k+1}= A_k z_k+B_k v_k,~ {z}_0 = \bar x_0,\label{eq:sls_nom} \\
& k = 0,\ldots, N-1\nonumber, \\
&\Px_{k+1,j} = A_k \Px_{k,j} + B_k \Pu_{k,j},~\Px_{j+1,j} = E_{j},\label{eq:SLP}\\
& j = 0,\ldots, N-1,~k = j+1,\ldots, N-1 \nonumber, \\
& \sum_{j=0}^{k-1}\left\|g_{k,i}\T \P_{k, j}\right\|_2+g_{k,i}^{\top}({z}_k, {v}_k)+b_{k,i} \leq 0,\label{eq:SLC}\\
&  i =1, \ldots, n_c,~ k=0,\ldots,N-1 \nonumber,\\
& \sum_{j=0}^{N-1}\left\|g_{\f,i}^{\top} \P_{N, j}\right\|_2+g_{\f,i}^{\top}{z}_N+b_{\f,i} \leq 0,\label{eq:SLCf}\\
& i =1, \ldots, \nf, \nonumber
\end{align}
\end{subequations}
with $\P_{k,j} \defeq (\Px_{k,j},\Pu_{k,j})$, $\P_{N,j} \defeq \Px_{N,j}$. Conditions \eqref{eq:sls_nom} implement nominal dynamics, while the tightened constraints~\eqref{eq:SLC}--\eqref{eq:SLCf} ensure the constraints~\eqref{eq:constraints} are satisfied robustly (see Proposition~\ref{prop:constraints}).
The cost~\eqref{eq:SLO} also minimizes the uncertainty propagation via the closed-loop response\footnote{This cost corresponds to the expected value of the variance of $J(\X,\U)$ under stochastic noise $w_k$, as detailed in~\citep{goulart2007affine}.}$\P$,
    \begin{equation}
    \begin{aligned}
        \label{eq:regulizer}
        &\tilde H_0(\P)\defeq \\
&\sum_{j=0}^{N-1} \left( \|\bar \Qf^{\frac{1}{2}} \Px_{N,j}\|_\F^2  + \sum_{k=j}^{N-1}  \left(\| \bar Q^{\frac{1}{2}} \Px_{k,j}\|_\F^2+ \| \bar R^{\frac{1}{2}} \Pu_{k,j}\|_\F^2 \right)\right),
    \end{aligned}
    \end{equation}
with $\bar Q\in \mathbb{S}_{++}^{\nx}$, $\bar R\in\mathbb{S}_{++}^{\nu}$, $\bar P\in \mathbb{S}_{++}^{\nx}$, and the high-dimensional variable $\P$ collecting all $\Px$, and $\Pu$. While this \gls{SLS} formulation is equivalent to disturbance feedback~\citep{goulart2006optimization}, its specificity is to consider explicitly the disturbance propagation~\eqref{eq:SLP}.

\begin{prop}
\label{prop:constraints}
\citep{goulart2006optimization}
System~\eqref{eq:LTV} in closed-loop with~\eqref{eq:disturbance_feedback} satisfies constraints~\eqref{eq:constraints} for all disturbances $w = (w_0, \dots, w_{N-1}) \in \mathcal{E}_{\nx} \times \cdots \times \mathcal{E}_{\nx}$ if and only if $(\Z,\V)$ and $(\Px,\Pu)$ satisfy the constraints in~\eqref{eq:sls}.
\end{prop}
\begin{remark}
Under suitable terminal ingredients $\Qf$, $g_{\f,i}$, $b_{\f,i}$, the optimization problem~\eqref{eq:sls}, used in a receding-horizon fashion, is recursively feasible and ensures a stable closed-loop system, see, e.g.,~\citep{goulart2006optimization}.
\end{remark}
\begin{remark}
In the \gls{SLS} literature, the constraint~\eqref{eq:SLP} is often written equivalently in its matrix form:
    \begin{equation}
\label{eq:slp}
\left[ \mathbf{I} - \mathbf{Z}{\mathbf{A}},~ - \mathbf{Z}{\mathbf{B}}  \right]  \begin{bmatrix}
     \Px\\
     \Pu\\
    \end{bmatrix} = \mathbf E,
\end{equation}
see, e.g.,~\citep[Eq.(2.7)]{anderson2019system}.
\end{remark}
\begin{remark}
The disturbance feedback~\eqref{eq:disturbance_feedback} can be equivalently reformulated as a state feedback, see, e.g., \citep{goulart2006optimization}.
\end{remark}
In the following, we exploit the structural similarity between the optimization problem~\eqref{eq:sls}, and the optimization problem efficiently solved in~\citep{messerer2021efficient}.
In contrast to standard tube-based \gls{MPC}, the SLS-based formulation~\eqref{eq:linear_SLS} is not conservative (Proposition~\ref{prop:constraints}). However, this approach results in a large number of decision variables, e.g., $\frac{N\cdot (N-1)}{2} \nx \nu$ for $\Pu$.
This paper proposes a structure-exploiting solver for the optimization problem~\eqref{eq:sls}.

\section{Efficient algorithm for SLS}
\label{sec:riccati_linear_sls}
In this section, we show how to efficiently solve the \gls{SLS}-based problem~\eqref{eq:sls}. 
The \gls{KKT} conditions~\citep{WrightStephenandNocedal1999NumericalOptimization} of~\eqref{eq:sls} are decomposed into two subsets.
The proposed algorithm iteratively solves each subset, corresponding to the nominal trajectory optimization and controller optimization employing efficient Riccati recursions (cf. Fig.\ref{fig:iterative_update_sls}).

\subsection{SLS and its KKT conditions}
\label{sec:KKT_SLS}
To write the optimization problem~\eqref{eq:sls} concisely, we introduce the auxiliary variables $\beta_{k,j}\in \R^{\nc}$ and the corresponding equality constraints
\begin{equation}
\label{eq:bo_def}
\beta_{k,j} = \tilde H_{k,j}(\P),
\end{equation}
similar for $\beta_N$, where we define
\begin{equation}
 \tilde H_{k,j}(\P)\defeq ( \|g_{k,1}\T \P_{k,j}\|_2^2, \ldots, \|g_{k,\nc}\T \P_{k,j}\|_2^2).
\label{eq:H_tilde_def}
\end{equation}

Based on these auxiliary variables $\beta_{k,j}$, the constraint tightening terms in~\eqref{eq:SLC}-\eqref{eq:SLCf} result as
\begin{equation}
   h\bo_k(\beta_k) \defeq \sum_{j=0}^{k-1} \sqrt{\beta_{k,j} + \epsilon_\beta},~h\bo_0\defeq 0_{\nc},\label{eq:def_hbok} 
\end{equation}
with
$\beta_k \defeq (\beta_{k,0}, \dots, \beta_{k,k-1})$. We note that the tightenings~\eqref{eq:def_hbok} are equivalent to those used in~\eqref{eq:sls} if $\epsilon_\beta=0$. However, in the following, we use a fixed $\epsilon_\beta>0$ to circumvent points of non-differentiability.

We summarize the \gls{SLS} problem~\eqref{eq:sls} as
\begin{subequations}
    \label{eq:linear_SLS}
    \begin{align}
    \min_{ y, \P,\beta} \quad &  J (y) + \tilde H_0(\P), \\
    \text{s.t.}\quad   & f(y) =0\label{eq:compact_sld}, \\
    & D(\P) =0, \label{eq:compact_slp}\\
    &  h(y) + h\bo(\beta) \le 0,\label{eq:compact_constraint}\\
& \tilde H(\P) -\beta =0, \label{eq:compact_bo}
    \end{align}
\end{subequations}
where $\beta \defeq (\beta_1, \dots, \beta_N)$, and $y=(\Z, \V) \in \mathbb{R}^{\ny}$ contains the variables associated with the nominal trajectory, the constraint~\eqref{eq:compact_sld} encodes the nominal dynamics~\eqref{eq:sls_nom}, and the constraint~\eqref{eq:compact_slp} encodes the disturbance propagation~\eqref{eq:SLP}. 
The tightened constraints~\eqref{eq:SLC}-\eqref{eq:SLCf} are split into two constraints: Inequality~\eqref{eq:compact_constraint} captures the nominal constraints $h(y)$ with a constraint tightening $h\bo(\beta) \defeq (0_{\nc},h\bo_1(\beta_1), \dots, h\bo_{N}(\beta_N))$ and Equality \eqref{eq:compact_bo} captures the nonlinear relation~\eqref{eq:bo_def}--\eqref{eq:H_tilde_def} between $\beta$ and the disturbance propagation $\P$.

To derive an efficient algorithm for solving the \gls{SLS}-based problem~\eqref{eq:linear_SLS}, we first condense the disturbance propagation~\eqref{eq:compact_slp} within the tightening definition~\eqref{eq:compact_bo}. 
This is achieved by recursively substituting all $\Px_{k,j}$ with a corresponding function of the gains $\Pu_{k,j}$ and disturbance gains $ E_j$ using~\eqref{eq:tube_propagation}. After substitution, we obtain the equivalent, more compact, optimization problem
\begin{subequations}
    \label{eq:linear_SLS_compact}
    \begin{align}
    \min_{ y, \M ,\beta} \quad &  J (y) + H_0(M), \\
    \text{s.t.}\quad   &   f(y) =0,\label{eq:linear_SLS_compact_lin_dyn} \\
    &  h(y) + h\bo(\beta) \le 0,\label{eq:linear_SLS_compact_bo}\\
    & H(\M) -\beta =0\label{eq:linear_SLS_compact_beta},
    \end{align}
\end{subequations}
where $\M$ collects all the disturbance feedback controller gains $\Pu_{k,j}$ in a vector, and $H$ is a composition function condensing the disturbance propagation~\eqref{eq:compact_slp}, such that \mbox{$H(\M) = \tilde H(\P)$} and $H_0(M) = \tilde H_0(\P)$, for any $\P$ with $D(\P) =0$.
The Lagrangian of~\eqref{eq:linear_SLS_compact} is given by
    \begin{equation}
    \begin{aligned}
            &\mathcal{L}(y, \beta, \M, \lambda, \mu, \eta) \\
            &= J (y) + H_0(M) + \lambda^\top f(y) + \mu^\top (h(y) + h\bo(\beta))\\
            &\quad  + \eta^\top (H(\M) - \beta),
    \end{aligned}
\end{equation}
with the dual variables $\lambda$, $\mu$, and $\eta$. The \gls{KKT} conditions of~\eqref{eq:linear_SLS_compact} are given by
\begin{subequations}
\label{eq:KKT}
\begin{align}
    \nabla J(y) + \nabla f(y) \lambda + \nabla h(y) \mu  &= 0, \label{eq:KKT_nab_y}\\
    \nabla h\bo(\beta) \mu - \eta &= 0,  \label{eq:KKT_nab_beta} \\
    \nabla H_0(M) + \nabla H(\M) \eta &= 0,\label{eq:KKT_nab_M}\\
     f(y) &= 0,\label{eq:prima_feas_KKT}\\
    0\leq \mu \perp h(y) + h\bo(\beta) &\leq 0,\label{eq:KKT_tight_constr} \\
    H(M) - \beta &= 0\label{eq:backoff_defs},
\end{align}
\end{subequations}
where $\perp$ denotes a complementarity condition. The next section shows how to exploit the structure in~\eqref{eq:KKT} to efficiently solve the SLS-based problem.
\subsection{Iterative trajectory and controller optimization}
In the following, we derive the proposed algorithm and we outline its convergence properties. The main idea of the proposed algorithm is to partition the~\gls{KKT} conditions~\eqref{eq:KKT} into two subsets solved alternatingly: one subset corresponding to the nominal trajectory optimization and the other to controller optimization.

In particular, the subset of necessary conditions~\eqref{eq:KKT_nab_y}, \eqref{eq:KKT_nab_beta}, \eqref{eq:prima_feas_KKT}, \eqref{eq:KKT_tight_constr} with fixed~$\bar \beta$, corresponds to a nominal trajectory optimization,
\begin{subequations}
\label{eq:KKT_nominal}
\begin{align}
    \nabla J (y) + \nabla f(y) \lambda + \nabla h(y) \mu &= 0, \\
    f(y) &= 0,\\
    0\leq \mu \perp h(y) + h\bo(\bar \beta) &\leq 0,\\
    \nabla h\bo(\bar \beta) \mu - \eta &= 0.~\label{eq:KKT_nominal_mu}
\end{align}
\end{subequations}
Conditions~\eqref{eq:KKT_nominal} yield an optimal dual $\bar \eta$. This optimal dual is fixed to solve the remaining necessary conditions~\eqref{eq:KKT_nab_M},~\eqref{eq:backoff_defs}, corresponding to the controller optimization
\begin{subequations}
\label{eq:KKT_contro}
\begin{align}
    \nabla H_0(M) + \nabla H(\M) \bar \eta &= 0,\label{eq:KKT_contro_eta}\\
    H(M) - \beta &= 0.\label{eq:KKT_contro_beta}
\end{align}
\end{subequations} 
The solution to~\eqref{eq:KKT_contro} yields a new update for $\beta$, used in~\eqref{eq:KKT_nominal}. 
A similar iterative scheme was derived in~\citep{messerer2021efficient} for a different robust optimal control problem, with each iteration refining the trajectory based on the computed feedback, while ensuring convergence towards an optimal solution (see Proposition~\ref{prop:conv_socp}).
Next, we detail how the decomposition in subsets~\eqref{eq:KKT_nominal} and~\eqref{eq:KKT_contro} enables a fast solution.
\subsection{Algorithm}

The solution of~\eqref{eq:KKT_nominal} can be obtained by solving a nominal trajectory optimization problem,
\begin{subequations}
    \label{eq:nominal_traj_opt_full}
    \begin{align}
    \min_{\Z, \V} \quad &J(\Z,\V), \\
    \text{s.t.}\quad  & {z}_{k+1}= A_k z_k+B_k v_k,~ {z}_0 = \bar x_0,\label{eq:sls_nom_full} \\
& h\bo_k (\bar \beta_k) + G_k({z}_k, {v}_k)+b_k \leq 0,\label{eq:SLC_full}\\
& k = 0,\ldots, N-1,\nonumber \\
& h\bo_N (\bar \beta_N)+G_{\f}{z}_N+b_{\f} \leq 0,\label{eq:SLCf_full}
    \end{align}
\end{subequations}
with the tightenings based on a fixed $\bar \beta$, and where $G_k$ and $G_\f$ collect, respectively, all $g_{k,i}\T$ and $g_{\f,i}\T$ in a matrix. 
Importantly, the stage-wise optimal control structure of the \gls{QP}~\eqref{eq:nominal_traj_opt_full} enables the use of tailored \gls{QP} solvers~\citep{kouzoupis2018recent} including efficient Riccati recursions~\citep{steinbach1995fast,rao1998application,frison2020hpipm}.
The optimal value of $\eta$ is then obtained by evaluating~\eqref{eq:KKT_nominal_mu}, i.e.,
\begin{equation}
    \eta_{k,j} = \frac{  \mu_k }{2\sqrt{ \bar \beta_{k,j}}}~j=0,\ldots,N-1,~k=j,\ldots,N-1,
   \label{eq:backoffs_duals}
\end{equation}
where $\eta_{k,j}\in \R^{n_c}$ are the dual variables that correspond to each $\beta_{k,j}$ for the equality constraint~\eqref{eq:backoff_defs}, and $\mu_k\in \R^{n_c}$ are the dual variables that correspond, at time step $k$, to the tightened constraints~\eqref{eq:KKT_tight_constr}. 

Then, to efficiently solve~\eqref{eq:KKT_contro} given a fixed $\bar \eta$, we notice that the stationarity condition~\eqref{eq:KKT_contro_eta} is linear in $M$ and has a unique solution. 
This condition can be solved with efficient Riccati recursions, as detailed in Section~\ref{sec:controller_Riccati}. Finally, the updated tightenings are obtained by evaluating~\eqref{eq:KKT_contro_beta}, which corresponds to the forward simulation given by~\eqref{eq:SLP}, and~\eqref{eq:bo_def}. Algorithm~\ref{algo} summarizes the steps previously outlined.
\begin{algorithm}
\caption{fast-SLS}\label{alg:cap}
\begin{algorithmic}
\Require $Q$, $R$, $P$, $A_k$, $B_k$, $E_k$, $G_{k}$, $b_{k}$ $G_\f$, $b_\f$
\State $\beta, \bar x_0  \gets \text{Initialize tightening (e.g., $\epsilon_\beta$), current state}$
\While{KKT~$\eqref{eq:KKT}$~\texttt{not satisfied}}
\State $\Z,\V,\mu \gets \text{Optimize nominal trajectory}$ \Comment{\eqref{eq:nominal_traj_opt_full}}
\State $\eta,C_{k,j} \gets \text{Update dual, and cost} $\Comment{\eqref{eq:backoffs_duals},~\eqref{eq:cost_update}}
\State $\Px,\Pu \gets\text{Optimize controller}$ \Comment{\eqref{eq:riccati}}
\State $\beta \gets~\text{Update tightening}$\Comment{\eqref{eq:SLP},~\eqref{eq:bo_def}}
\EndWhile\\
\textbf{Output:} $\Z^\star,~\V^\star,~{\Px}^\star,~{\Pu}^\star$ \text{satisfy~\eqref{eq:linear_SLS}}
\end{algorithmic}
\label{algo}
\end{algorithm}

\begin{prop}
\label{prop:conv_socp}
Denote by $s^\star \defeq ( y^\star, M^\star, \beta^\star )$ the unique optimal solution of~\eqref{eq:linear_SLS} or~\eqref{eq:linear_SLS_compact} and assume standard regularity conditions hold at $s^\star$, i.e., linear independence constraint qualification (LICQ), the second order sufficient condition (SOSC), and strict complementarity. 
    Then, the sequences $\{ s_i\}_{i=0}^{\infty}$, generated by Algorithm~\ref{alg:cap}, converge to $s^\star$ if $\sigma = \max_k\|E_k\|_2$ is sufficiently small and if $\beta$ and the corresponding $y$ are initialized sufficiently close to $s^\star$.
    The convergence rate is q-linear in the nominal trajectory variables $y$ and r-linear in the controller variables $\P$, {and $\beta$}.
\end{prop}
\begin{pf}
The proof is analogous to~\citep[Thm. 10]{messerer2021efficient}. 
Specifically, by fixing $\bar\beta$ in~\eqref{eq:KKT_nominal}, we implicitly neglect the corresponding Jacobians with respect to the remaining variables. As in~\citep{messerer2021efficient}, the magnitude of the neglected Jacobians is proportional to $\sigma$. Hence, for $\sigma$ small enough, we obtain linear convergence, whereas quadratic convergence
is not possible because of the neglected gradients.
Notably, \cite{messerer2021efficient} use a gradient-correction to ensure convergence to a KKT point.
However, this is not needed for the convex problem~\eqref{eq:linear_SLS} and Algorithm 1 converges to the global minimizer $s^\star$.
$\hfill\square$  
\end{pf}
\begin{remark}
    Without inequality constraints, a single trajectory and Riccati recursion for the controller optimization yields an optimal solution.
\end{remark}

\begin{remark}
    As stopping criterion, we compare the norm of the difference between two consecutive primal solutions to a specified threshold, denoted as $\epsilon_m$.
\end{remark}

\begin{remark}
The solution of~\eqref{eq:nominal_traj_opt_full} is feasible with respect to the \gls{SOCP}~\eqref{eq:sls} at every iteration of the iterative scheme Algorithm~\ref{algo}, allowing for safe early termination.
\end{remark}

Next, we focus on exploiting the structure of~\eqref{eq:KKT_contro_eta}.

\subsection{Controller optimization via efficient Riccati recursions}
\label{sec:controller_Riccati}
We exploit the numerical structure in~\eqref{eq:KKT_contro_eta}, to efficiently compute a solution. Equation~\eqref{eq:KKT_contro_eta} corresponds to first-order necessary conditions of
\begin{subequations}
    \label{eq:tube_update}
    \begin{align}
    \min_{ \P} \quad &  \eta\T \tilde H(\P)+\tilde H_0(\P), \label{eq:tube_cost}\\
    \text{s.t.}\quad   & D(\P) =0.
    \end{align}
\end{subequations}

The structure of the cost function~\eqref{eq:tube_cost} can be interpreted intuitively: the tightenings $\|g_{k,i} \T \P_{k,j}\|_{2}^2$ associated with larger dual variables $\eta$ in~\eqref{eq:KKT_nab_M} are penalized more heavily, leading to stronger disturbance rejection in that direction.

After some algebraic manipulations, the optimization problem~\eqref{eq:tube_update} can be written equivalently as
\begin{subequations}
    \label{eq:tube_update_full}
    \begin{align}
    \min_{{\Px}, {\Pu}} \quad & \sum_{j=0}^{N-1}   \| C_{N,j} \Px_{N,j}\|_{\mathcal{F}}^2  + \sum_{k=j}^{N-1} \|  C_{k,j} \P_{k,j}\|_{\mathcal{F}}^2, \label{eq:tube_update_full_cost}\\
    \text{s.t.}  ~&\Px_{k+1,j} = A_{k} \Px_{k,j} +B_{k} \Pu_{k,j},\\
    &{\Px_{j+1,j} = E_j},\\
    & j = 0,\ldots, N-1,~k = j+1,\ldots, N-1,\nonumber 
    \end{align}
\end{subequations}
where we define the concatenation
{\begin{equation}
\begin{aligned}
    &{C}_{k,j} \defeq {\left( \text{diag}\left(\sqrt{ \eta_{k,j}}\right) G_k, \begin{bmatrix}
       \bar Q^{1/2}& \mathbf{0}_{\nx,\nu} \\
    \mathbf{0}_{\nu,\nx} & \bar R^{1/2}
    \end{bmatrix}\right)},\\
    &     {C}_{N,j} \defeq{ \left( \text{diag} \left(\sqrt{ \eta_{N,j}}\right) G_\f , \bar P^{1/2}\right)}.
\label{eq:cost_update}
\end{aligned}
\end{equation}
Each term of $\tilde H_0(\P)$ in~\eqref{eq:tube_cost} has been incorporated in the sum~\eqref{eq:tube_update_full_cost} by leveraging Frobenius norm properties.

We introduce the block decomposition
\begin{equation}
\begin{bmatrix}
        C^\x_{k,j} & C^{\x\u}_{k,j}\\
        C^{\x\u}_{k,j} & C^\u_{k,j}\\
    \end{bmatrix} \defeq   C_{k,j} \T C_{k,j},
\end{equation}

such that $C^\x_{k,j} \in \mathbb{S}_{++}^{\nx}, C^\u_{k,j} \in \mathbb{S}_{++}^{\nu}$ and $C^{\u\x}_{k,j}={C^{\x\u}_{k,j}}\T$ of corresponding dimension, for $j=0,\ldots, N-1,~ k=j, \ldots, N-1$, and similar for $C^\x_{N,j}$. The following Theorem shows that the solution to~\eqref{eq:tube_update_full} is given by $N$ independent Riccati recursions.
\begin{thm}
\label{prop:riccati}
Problem~\eqref{eq:tube_update_full} has a unique minimizer which is given by $N$ parallel  backward Riccati recursions,
    \begin{equation}
\begin{aligned}
\label{eq:riccati}
S_{N,j} & =C^\x_{N,j}, \\
K_{k,j} & =-\left(C^\u_{k,j}+B_k^{\top} S_{k+1,j}B_k\right)^{-1}\left(C^{\u\x}_{k,j}+B_k^{\top} S_{k+1,j}A_k\right), \\
S_{k,j} & = C^\x_{k,j}+A_k^{\top} S_{k+1,j}A_k+\left(C^{\x\u}_{k,j}+A_k^{\top} S_{k+1,j}B_k\right) K_{k,j},
\end{aligned}
    \end{equation}
followed by $N$ parallel forward propagations,
\begin{equation}
\label{eq:lyap_forward}
\begin{aligned}
\Px_{j+1,j}&= E_j\\
    \Pu_{k,j }&= K_{k,j} \Px_{k,j},~ \Px_{k+1,j}= (A_k +B_k K_{k,j})\Px_{k,j},
\end{aligned}
    \end{equation}
for $j = 0,\ldots, N-1,~ k=j+1, \ldots, N-1$.
\end{thm}

\begin{pf}
First, we note that each disturbance feedback controller operates independently, allowing us to address the solution for each disturbance $ j=0,\ldots, N-1$ as follows:
\begin{subequations}
    \label{eq:single_tube_update}
    \begin{align}
    \min_{\substack{\Px_{j+1:N, j},\\ \Pu_{j+1:N, j}}}~& \| C_{N,j}\Px_{N,j}\|_\F^2  + \sum_{k=j+1}^{N-1}  \| C_{k,j} \P_{k,j}\|_\F^2, \\
    \text{s.t.}  ~ &\Px_{j+1,j} = E_j,\\
    &\Px_{k+1,j} = A_{k} \Px_{k,j} +B_{k} \Pu_{k,j},\\
    &~ k = j+1,\ldots, N-1.\nonumber
    \end{align}
\end{subequations}
We have a quadratic cost,
\begin{equation}
    \begin{aligned}
        \| C_{k,j} \P_{k,j}\|_\F^2 &= \text{Trace}\left( {\P_{k,j}}\T C_{k,j}\T C_{k,j} {\P_{k,j}} \right)\\
        &= \sum_{i=1}^{\nx} e_i\T {\P_{k,j}}\T C_{k,j}\T C_{k,j} \P_{k,j} e_i,
    \end{aligned}
\end{equation}
with $e_i\T$ is the $i^\text{th}$ row of the identity matrix $I_{\nx}$, and a linear unconstrained dynamics, as in~\citep{tseng2020system}, with matrix dynamics. Notably, each column of $\Px_{k,j}$, $\Pu_{k,j}$, i.e., $\P_{k,j}e_i$, can be interpreted as independent states and input for each $i=1,\ldots,\nx$, with identical cost and dynamics but different initial conditions.
Hence, Problem~\eqref{eq:single_tube_update} corresponds to $\nx$ \gls{LQR} problems with different initial conditions and hence can be efficiently solved via a single identical Riccati recursion~\eqref{eq:riccati} and a forward matrix propagation~\eqref{eq:lyap_forward}.
$\hfill\square$\end{pf}
\begin{remark}
\label{rem:extension_matrix_E}
For disturbances $w_k\in \R^{\nw}$ with $\nw < \nx$, Theorem~\ref{prop:riccati} can be extended to left-invertible $E_k \in \R ^{\nx  \times n_\mathrm{w}}$ as used in\citep{herold2022computationally}\footnote{A left-invertible matrix $E_k$ implies there exists a matrix $E_k^\dagger \in \mathbb{R}^{n_\mathrm{w} \times \nx}$ such that $E_k^\dagger E_k = I_{\nw}$.}. In this case,  Problem~\eqref{eq:tube_update_full} is optimally solved by the Riccati recursions~\eqref{eq:riccati} and replacing~\eqref{eq:lyap_forward} by $N$ parallel forward propagations
\begin{equation}
\label{eq:lyap_forward_dagger_E}
\begin{aligned}
\Px_{j+1,j}&= E_j, ~ \Omega_j = E_j E_j^\dagger,\\
    \Pu_{k,j }&= K_{k,j} \Omega_j \Px_{k,j},~ \Px_{k+1,j}= (A_k +B_k K_{k,j})\Px_{k,j},
\end{aligned}
    \end{equation}
for $j = 0,\ldots, N-1,~ k=j+1, \ldots, N-1$. Note that for $E_j$ invertible,~\eqref{eq:lyap_forward_dagger_E} reduces to~\eqref{eq:lyap_forward} with $\Omega_j = I_{\nx}$.
\end{remark}
\begin{remark}
\label{rmq:nonlinear_extension}
The proposed algorithm can be extended to address a nonlinear robust MPC formulation~\citep{Leeman2023RobustSynthesis}. To efficiently solve the resulting {NLP}, a modified {SCP} can be utilized, solving a series of convex problems iteratively to approximate the original nonconvex problem. Each iterate in the sequence is equivalent to a convex \gls{SLS}-based problem~\eqref{eq:sls}, which can be efficiently addressed using the proposed algorithm. Crucially, all the constraints of the convexified problems have the same numerical structure as in the linear case, and can hence be handled using Algorithm 1. 
\extension{
A more detailed derivation can be found in Appendix~\ref{sec:nonlinear_extension}.
}{
A more detailed derivation can be found in an extended version of this paper~\cite[App.~A]{leeman2024}.
}
\end{remark}
\begin{figure*}[ht!]
    \centering
    \includegraphics[width = 0.9\textwidth]{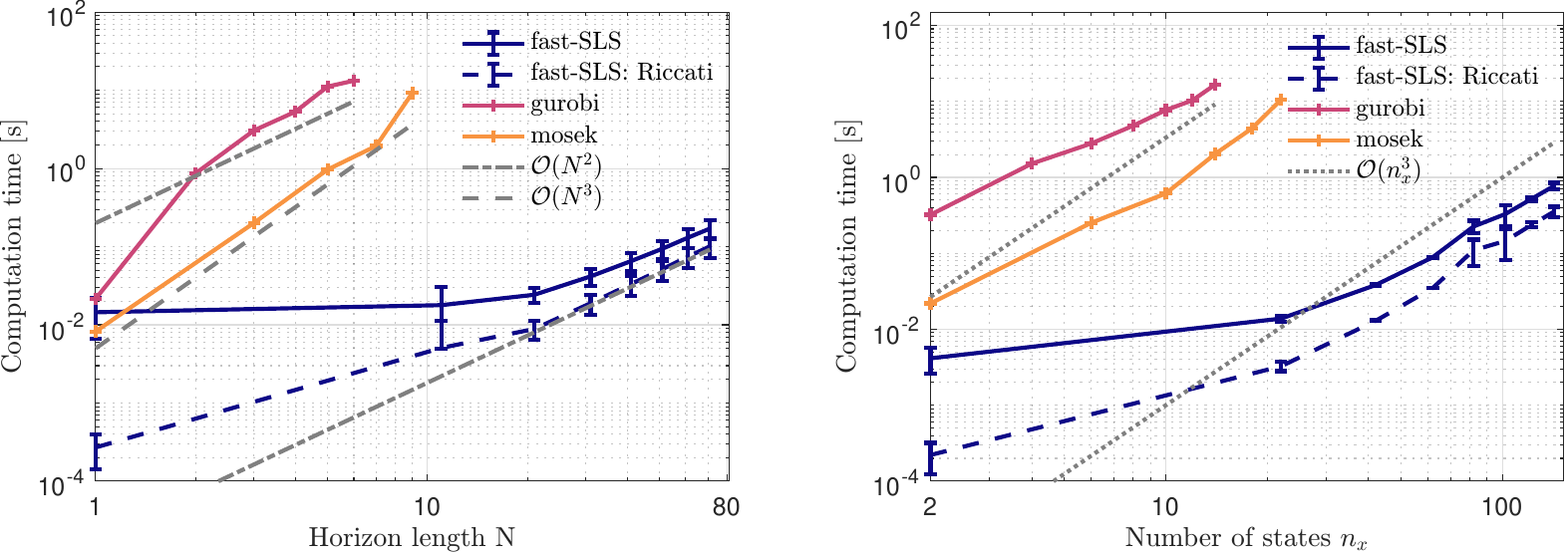}
    \caption{Scalability of different solvers with horizon length for $L=10$ masses (left). Scalability of different solvers with the number of states for horizon $N=10$ (right).
    The proposed algorithm's (fast-SLS) computation times, and the Riccati recursions contributions~\eqref{eq:riccati} (fast-SLS: Riccati) are shown with standard deviations.}
    \label{fig:fig4}
\end{figure*}
\begin{figure}[ht!]
    \centering
    \includegraphics[width=0.45\textwidth]{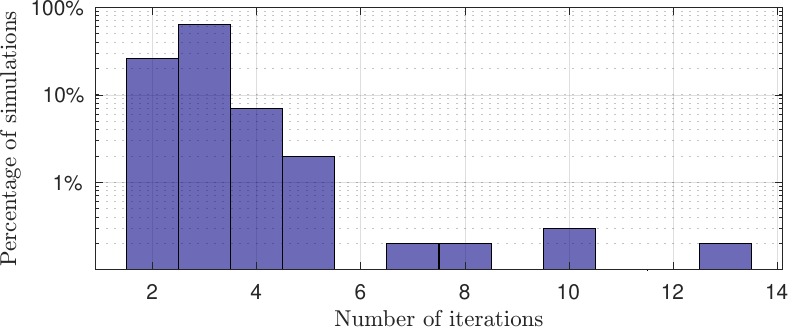}
    \caption{Number of iterations required until convergence for $10^3$ randomly sampled initial conditions, for horizon $N=25$, and $L=25$ masses.}
    \label{fig:fig3}
\end{figure}
\subsection{Scalability}
\label{sec:ccalability}
For unconstrained linear quadratic problems, the Riccati recursion scales as $\mathcal{O}(N(\nx^3+\nu^3))$~\citep{rawlings2017model}.
Using a naive alternative based on state augmentation to solve equation~\eqref{eq:sls}, the number of states in each stage would grow from $\nx$ to $\nx + N\nx^2$, with a similar expansion for the inputs. Consequently, this naive implementation approach leads to a computational scalability of $\mathcal{O}(N^4 ( \nx^3 + \nu^3 )\nx^3 )$. The solver introduced by \cite{goulart2008efficient} leverages the sparsity structure for factorization of the KKT matrix, leading to an improved scaling of $\mathcal{O}(N^3(\nx^3 + \nu^3))$.

In contrast, the proposed algorithm alternatively solves two subsets of the KKT conditions. Hence, the computation times are dominated by the $N$ Riccati recursions executed at each iteration. If implemented sequentially, this scales as $\mathcal{O}(N^2(\nx^3+\nu^3))$. 
By parallelizing the Riccati recursions across $N$ cores, the scalability can be reduced to $\mathcal{O}(N(\nx^3+\nu^3))$, achieving the same computational scalability as the nominal trajectory optimization.

\section{Numerical example}

\label{sec:numerical}
In this section, we demonstrate the computational speedup achieved by the proposed algorithm compared to available alternatives.
We show how the proposed algorithm scales with increasing horizon and increasing number of states. The nominal trajectory optimization~\eqref{eq:nominal_traj_opt_full} is solved with \texttt{osqp}~\citep{osqp}, via its \texttt{CasADi} interface~\citep{Andersson2019}, while the Riccati recursions~\eqref{eq:riccati} are implemented using \texttt{MATLAB}\footnote{An open-source implementation is available at\\ \url{https://gitlab.ethz.ch/ics/fast-sls},\\ doi: \url{https://doi.org/10.3929/ethz-b-000666883}}.

The proposed algorithm is compared against directly solving the \gls{SOCP}~\eqref{eq:sls} using general-purpose commercial solvers,~\cite{gurobi} and Mosek~\cite{mosek}, both interfaced using \texttt{YALMIP}~\citep{Lofberg2004}. Computations are carried out on a Macbook Pro equipped with M1 processor with 8 cores and 16 GB of RAM, running macOS Sonoma.

In the following example, we consider a chain of $L$ mass-spring-damper systems as in~\citep{domahidi2012efficient}, corresponding to a linear system with $\nx = 2 \cdot L$, fixed at one end, with mass $m=1$, spring constant $k=10$, and damper constant $d=2$, and a force actuation on each mass. The system is discretized using a time step $\Delta t = 0.1$. We use $Q = \Qf=3 {I_{\nx}}$, $R = { I_{\nu}}$, $E= 0.1 I_{\nx}$. 
The states $x$ and the force input $u$ are constrained at each time step as $\|x\|_\infty \le 4$, and $\| u\|_\infty \le 4$.
We use a convergence tolerance on the primal variables $y$ of $\epsilon_m=10^{-8}$, and a minimal tightening of $\epsilon_\beta = 10^{-10} \bm{1}_{\nc}$.
The computation times of the proposed algorithm are calculated based on an average of $30$ randomly sampled feasible initial conditions. We also show the computation times for the Riccati recursions alone~\eqref{eq:riccati}.

\subsubsection{Scalability with state dimension}
For horizon $N=10$, and an increasing number of masses $L$, we investigate the computation times for different solvers. 
In Fig.~\ref{fig:fig4}~(right), we see that the proposed algorithm improves computation times compared to Mosek and Gurobi, especially for systems with high state dimension, achieving an improvement in computation speed by a factor of up to $10^3$ times. However, as predicted, our method scales as $\mathcal{O}(\nx^3)$, which is the same asymptotic scaling as for Mosek and Gurobi.

\subsubsection{Scalability with horizon length}
We investigate the computation times for different solvers, using $L=10$ masses and increasing horizon length $N$. In Fig.~\ref{fig:fig4} (left), {we see that our method scales as $\mathcal{O}(N^2)$}. In comparison, Mosek and Gurobi scale respectively approximately as  $\mathcal{O}(N^3)$, and $\mathcal{O}(N^2)$, {indicating the use of a sparse factorization method, as in~\citep{goulart2008efficient}}.
Crucially, the proposed algorithm improves computation times, especially for long horizons, achieving an improvement in computation speed by a factor of up to $10^3$ times.

\subsubsection{Convergence}
For horizon $N=25$, and $L=25$ masses, we evaluate the number of iterations required to numerically converge for $10^3$ randomly sampled feasible initial conditions for the proposed solver. In Fig.~\ref{fig:fig3}, we see that, in the vast majority of cases, only a few iterations are required to converge to the optimal solution.
All the solvers considered converge to the same optimal solution.
\begin{remark}
As seen in Fig.~\ref{fig:fig4}, the computation times are dominated by the Riccati recursions~\eqref{eq:riccati}. It can be noted that due to overhead in the \texttt{MATLAB} implementation, the computation times are less reliable for small problems.
The computation speed of the proposed algorithm could benefit from a more efficient implementation, e.g., relying on the Cholesky decomposition of the cost matrices~\eqref{eq:cost_update}, compare~\citep{frison2016algorithms}.
\end{remark}

\section{Conclusion}

SLS is a flexible technique to jointly optimize the nominal trajectory and controller in robust MPC by using a disturbance feedback parametrization. 
The resulting SLS-based MPC yields robust constraint satisfaction and provides stability. In this paper, we address the main bottleneck of \gls{SLS}: the computation times. 
The proposed algorithm iteratively optimizes the nominal trajectory and then the controller. The nominal trajectory optimization is solved via an efficient QP solver and controller optimizations are solved via Riccati recursions.
Each iterate scales with $\mathcal{O}(N^2(\nx^3 + \nu^3))$. 
Leveraging results from~\cite{messerer2021efficient}, the proposed algorithm locally converges q-linearly to an optimal solution. 
In a numerical comparison, the proposed algorithm offers speed ups of up to $10^3$ compared to general-purpose commercial solvers. 
We expect that the proposed algorithm can be adapted to address variations of SLS, particularly those involving sparsity, structural constraints~\citep{tseng2020system} or polytopic disturbances and parameters~\citep{chen2022robust}.

\bibliography{references}     

\extension{
\newpage
\appendix
\section{Nonlinear Extension}
\label{sec:nonlinear_extension}

In the following, we show how the efficient solution to the SOCP~\eqref{eq:sls} can be extended to address nonlinear systems. 
Based on~\citep{Leeman2023RobustSynthesis,leeman2023_CDC}, we address robust control of nonlinear systems using the following nonlinear program
\begin{subequations}
\label{eq:nonlinear_sls}
\begin{align}
\min_{\substack{\mathbf{\Phi}_{\mathbf{x}},\mathbf{\Phi}_{\mathrm{u}},\\\mathbf{z}, \mathbf{v},\bm{\tau}}} &  J(\Z,\V) +  \tilde H_0(\P),\\
\text { s.t. } &  {z}_{k+1}= \phi(z_k,v_k),~ {z}_0 = \bar x_0, ~ k = 0,\ldots, N-1,\label{eq:nonlinear_nominal}\\
&\Px_{k+1,j} =\nabla_z \phi\left(z_k, v_k\right)\T  \Px_{k,j} + \nabla_v \phi\left(z_k, v_k\right)\T \Pu_{k,j},\label{eq:SLP_nonlinear}\\
&\Px_{j+1,j} = \alpha_1^{-1} E_j + {\alpha_2^{-1} \sqrt{\nx} \mu \tau_j^2},\label{eq:SLP_ic_nonlinear}\\
& j = 0,\ldots, N-1,~k = j+1,\ldots, N-1, \nonumber \\
& \sum_{j=0}^{k-1}\left\|e_i^{\top} \P_{k, j}\right\|_2-\tau_k \leq 0,\label{eq:SLC_lin}\\
&  i =1, \ldots, \nx + \nu, k= 0,\ldots,N -1, \nonumber\\
& \sum_{j=0}^{k-1}\left\|g_{k,i}^{\top} \P_{k, j}\right\|_2+g_{k,i}^{\top}\left({z}_k, {v}_k\right)+b_{k,i} \leq 0,\label{eq:SLC_nonlinear}\\
&  i =1, \ldots, n_c,~k= 0,\ldots,N -1, \nonumber\\
& \sum_{j=0}^{N-1}\left\|g_{\f,i}^{\top} \Px_{N, j}\right\|_2+g_{\f,i}^{\top}{z}_N+b_{\f,i} \leq 0,\label{eq:SLCf_nl}\\
&  i =1, \ldots, \nf,\nonumber
\end{align}
\end{subequations}
where $\bm{\tau} \in \R^{N}$ and $e_i\T$ is the $i^\text{th}$ row of the identity matrix $I_{\nx + \nu}$.
The function $\phi:\R^{\nx}\times \R^{\nu} \rightarrow \R^{\nx}$ is assumed to be three times continuously differentiable and is used to define the dynamics~\eqref{eq:nonlinear_nominal}, and the disturbance propagation~\eqref{eq:SLP_nonlinear}.
In~\eqref{eq:SLP_ic_nonlinear} the effect of the additive disturbance $w_k$ and the linearization error based on the worst-case curvature $\mu\in \R^{\nx \times \nx} $ of $\phi$ is lumped, see~\citep{Leeman2023RobustSynthesis}. Here, $\alpha_1, \alpha_2 \ge 0$ is used to over-approximate the sum of ellipsoids~\citep[Thm 2.4]{houska2011robust} with $\alpha_1 + \alpha_2 \leq 1$.
The constraints~\eqref{eq:SLC_lin} are used to optimize the dynamic linearization error over-bound $\tau_k$, while still jointly optimizing the nonlinear nominal trajectory, and feedback controller, such that robust constraint satisfaction is guaranteed for the nonlinear dynamics using~\eqref{eq:SLC_nonlinear}, and~\eqref{eq:SLCf_nl}.

We write the optimization problem~\eqref{eq:nonlinear_sls} compactly as
\begin{subequations}
    \label{eq:nonlinear_sls_compact}
    \begin{align}
    \min_{  y,\P,\beta} \quad &  {J( y) + \tilde H_0(\P) }, \\
    \text{s.t.}\quad   & f_\textrm{nl}( y)=0,\label{eq:compact_nonlinear}\\
    & D_\textrm{nl}( y,\P) =0, \label{eq:slp_nlp}\\
    & \tilde{h}( y) + h\bo(\beta) \le 0,\label{eq:nonlinear_sls_compact_beta_def}\\
& \tilde H(\P) -\beta =0, \label{eq:nonlinear_sls_compact_beta}
    \end{align}
\end{subequations}
where $ y\defeq (\Z,\V,\bm{\tau})$ collects all the nominal trajectory variables, including $\bm{\tau}$. The constraint~\eqref{eq:compact_nonlinear} corresponds to the nominal trajectory~\eqref{eq:nonlinear_nominal}.
The disturbance propagation $D_\textrm{nl}( y,\P)$ in~\eqref{eq:slp_nlp} is linear in $\P$, and \eqref{eq:nonlinear_sls_compact_beta_def},~\eqref{eq:nonlinear_sls_compact_beta} lump all the inequality constraints~\eqref{eq:SLC_lin},~\eqref{eq:SLC_nonlinear}, and~\eqref{eq:SLCf_nl}.

\begin{algorithm}[ht]
\caption{SCP-SLS}\label{alg:scp-sls}
\begin{algorithmic}
\Require $Q$, $R$, $P$, $G_{k}$, $b_{k}$, $G_\f$, $b_\f$
\State $\bar M, \bar \eta, \bar y  \gets \text{Initialize controller, duals, nominal trajectory}$
\While{\texttt{not converged}}
\State $\nabla f_\mathrm{nl}(\bar y), \Gamma \bar \eta \gets \text{Evaluate gradients}$ \Comment{\eqref{eq:reduced_linearized_SLS_a},~\eqref{eq:reduced_linearized_SLS_b}}
\State $\Delta y, \bar M \gets \text{Solve~\eqref{eq:reduced_linearized_SLS}} $\Comment{Alg.~\ref{algo}}
\State $\bar y \gets ~\text{Update nominal}$\Comment{$\bar y + \Delta y$ }
\EndWhile\\
\textbf{Output:} $\Z^\star,~\V^\star,~\bm{\tau}^\star,~{\Px}^\star,~{\Pu}^\star$ \text{satisfy~\eqref{eq:nonlinear_sls}}
\end{algorithmic}
\label{algo_scp}
\end{algorithm}

\newpage
After substituting~\eqref{eq:slp_nlp} in~\eqref{eq:nonlinear_sls_compact_beta}, similarly as in the linear case Section~\ref{sec:KKT_SLS}, we obtain an equivalent optimization problem
\begin{subequations}
    \label{eq:nonlinear_sls_compact}
    \begin{align}
    \min_{  y,M,\beta} \quad &  {J( y)} +  H_0(M), \label{eq:nonlinear_sls_compact_cost}\\
    \text{s.t.}\quad   & f_{\mathrm{nl}}( y)=0,\label{eq:nonlinear_sls_compact_dynamics}\\
    & \tilde h( y) + h\bo(\beta) \le 0,\label{eq:nonlinear_SLS_compact_constraint}\\
& H(M,y) -\beta \le 0 \label{eq:nonlinear_sls_compact_beta_2}.
    \end{align}
\end{subequations}

In~\eqref{eq:nonlinear_sls_compact_beta_2}, we replaced the equality constraint $ H(M, y)-\beta\leq 0$ by an \textit{inequality}, which is convex in $M$. {Likewise, the inequality~\eqref{eq:nonlinear_SLS_compact_constraint} can be losslessly convexified using auxiliary variables.} These modifications do not alter the optimal solution $y$ and $M$ of~\eqref{eq:nonlinear_sls_compact}, as the tightenings $\beta$ do not appear in the cost~\eqref{eq:nonlinear_sls_compact_cost}.
Applying a standard \gls{SCP} approach~\citep{dinh2010local}, we first linearize the nonconvex constraints in~\eqref{eq:nonlinear_sls_compact} at $\bar y$, $\bar M$, yielding
\begin{subequations}
    \label{eq:linearized_SLS}
    \begin{align}
    \min_{ \Delta y,  M,\beta} \quad &  J(\bar y+ \Delta y) + H_0(M) \label{eq:linearized_SLS_cost}\\
    \text{s.t.}\quad   & f_\mathrm{nl}(\bar y) + \nabla f_\mathrm{nl}\T(\bar y) \Delta y=0,\\
    & \tilde h(\bar y + \Delta y) + h\bo(\beta) \le 0,\label{eq:linearized_SLS_cons}\\
& H(M,\bar y) + \Gamma \Delta y -\beta \leq 0\label{eq:linearized_SLS_backoffs},
    \end{align}
\end{subequations}
with the gradient 
\begin{equation}
\Gamma \defeq \nabla_y H( \bar { M},\bar y)\T.
\end{equation}
\begin{remark}
In the SOCP~\eqref{eq:linearized_SLS}, the constraints~\eqref{eq:linearized_SLS_backoffs} only include the linearization with respect to $y$, preserving the advantageous numerical structure of the optimization problem.
\end{remark}

In the optimization problem~\eqref{eq:linearized_SLS}, the term $\Gamma \Delta y$ in~\eqref{eq:linearized_SLS_backoffs} prevents the direct use of Algorithm 1 to solve~\eqref{eq:linearized_SLS}. Hence, to be able to use Algorithm 1, and leverage the structure we use a variation of \gls{SCP}~\citep{bock2007constrained}, where the constraint~\eqref{eq:linearized_SLS_backoffs} is approximated with
\begin{equation}
    H( M, \bar y) -\beta \le 0.
\end{equation}
Due to this inexact Jacobian approximation~\citep{bock2007constrained}, a correction in the cost function~\eqref{eq:linearized_SLS_cost} is needed to maintain convergence properties of the SCP scheme.
In particular, we iteratively solve a series of optimization problems, equivalent to the SOCP~\eqref{eq:sls}:
\begin{subequations}
\label{eq:reduced_linearized_SLS}
    \begin{align}
    \min_{ \Delta y,  M,\beta} \quad &  J(\bar y + \Delta y)+ H_0(M) + \Gamma \bar \eta \Delta y , \label{eq:reduced_linearized_SLS_a}\\
    \text{s.t.}\quad   & f_\mathrm{nl}(\bar y) + \nabla f_\mathrm{nl}(\bar y)\T \Delta y=0,\label{eq:reduced_linearized_SLS_b}\\
    & \tilde h(\bar y + \Delta y) + h\bo(\beta) \le 0,\\
    &H(M, \bar y) -\beta \le 0,\label{eq:reduced_linearized_SLS_bo}
    \end{align}
\end{subequations}
and iteratively update the optimal nominal trajectory $\bar y$ and the optimal dual variable $\bar \eta$ associated to~\eqref{eq:reduced_linearized_SLS}.
Using this approximation, each \gls{SOCP} iterate~\eqref{eq:reduced_linearized_SLS} has the same structure as~\eqref{eq:linear_SLS_compact}, and, hence, can be efficiently solved using the proposed algorithm, described in Section~\ref{sec:riccati_linear_sls}. Indeed, the inequality in~\eqref{eq:reduced_linearized_SLS_bo} can be replaced by an equality without changing the optimal solution in $\Delta y$ and $M$. Algorithm~\ref{algo_scp} summarizes the steps previously outlined.

Under some standard regularity conditions, the \gls{SCP} based on~\eqref{eq:reduced_linearized_SLS} converges to a KKT point~\citep{messerer2021survey} of the original optimization problem~\eqref{eq:nonlinear_sls}.
In particular, the optimization problems~\eqref{eq:reduced_linearized_SLS} and~\eqref{eq:linearized_SLS} have the same stationarity conditions.
}{
}

\end{document}